\documentclass[12pt]{article}
\usepackage{amsmath,amsthm,amsfonts,amssymb,mathrsfs,bm,graphicx,amscd,epsfig,psfrag,verbatim,hyperref}
\usepackage[usenames]{color}
\usepackage{fancyhdr}
\parskip 	\bigskipamount

\newtheorem{theorem}{Theorem}[section]
\newtheorem{lemma}[theorem]{Lemma}
\newtheorem{corollary}[theorem]{Corollary}

\newtheorem{remark}{Remark}[section]

\newtheorem{definition}{Definition}
\newtheorem{claim}{Claim}

\def\Z{\mathbb{Z}}

\def\R{\mathbb{R}}
\def\C{\mathbb{C}}
\def\E{\mathbb{E}}
\def\P{\mathbb{P}}

\def\del{\delta}

\def\F{\mathcal{F}}

\def\La{\Delta}

\def\D{\mathcal{D}}
\def\eps{\epsilon}
\def\A{\mathcal{A}}
\def\o{\omega}

\def\z{\zeta}
\def\uz{\underline{\zeta}}
\def\us{\underline{s}}

\def\el{\mathcal{L}}
\def\ph{\varphi}
\def\S{\mathcal{S}}

\def\nat{\mathbb{N} \cup \{0\}}
\def\uz{\underline{\zeta}}
\def\uout{\Upsilon_{\out}}

\def\out{{\mathrm{out}}}

\def\els{\mathcal{L}_{\Sigma}}
\def\ol{\overline}

\def\c{\complement}
\def\M{\mathbb{M}}

\def\G{\mathcal{G}}
\def\F{\mathcal{Z}}

\def\nat{\mathbb{N} \cup \{0\}}
\def\uz{\underline{\z}}

\def\uout{\Upsilon_{\text{out}}}

\def\el{\mathcal{L}}
\def\els{\el_{\Sigma}}

\def\inn{\mathrm{in}}

\def\mm{\mathfrak{m}}

\def\ol{\overline}

\def\eps{\varepsilon}

\def\A{\mathcal{A}}

\def \del{\delta}
\def \z{\zeta}
\def \o{\omega}
\def \Z{\boldsymbol{Z}}

\def \d{\mathrm{d}}

\def \M{\boldsymbol{M}}

\def\E{\mathbb{E}}

\def\B{\mathcal{B}}
\def\zz{\underline{\zeta}}

\def\ul{\underline}
\def\La{\Lambda}
\def\uz{\underline{\zeta}}
\def\uza{{\uz_1}}
\def\uzb{{\uz_2}}
\def\Supp{\text{Supp}}
\def\m{\underline{m}}
\def\a{\alpha}
\def\aa{\underline{a}}
\def\M{\mathcal{M}}
\def\mm{\mathfrak{M}}

\renewcommand{\l}[0]{\left }
\renewcommand{\r}[0]{\right}

\makeatletter
\renewcommand*{\@cite@ofmt}{\hbox}
\makeatother

\begin{document}
\title{\bf Palm measures and rigidity phenomena in point processes }

\author{\href{http://www.princeton.edu/~sg18/}{Subhroshekhar Ghosh}\\ Princeton University\\ 
}

\maketitle

\begin{abstract}
 We study the mutual regularity properties of Palm measures of point processes, and establish that  a key determining factor for these properties is the rigidity-tolerance behaviour of the point process in question (for those processes that exhibit such behaviour). Thereby, we extend the results of \cite{OsSh}, \cite{Bu-3}, \cite{Ol} to new ensembles, particularly those that are devoid  of any determinantal structure. These include the zeroes of the standard planar Gaussian analytic function and several others. 
\end{abstract}

\section{Introduction}

Our aim in this article is to study the mutual singularity (and continuity) properties of Palm measures of point processes. Roughly speaking, the Palm measure of a point process $\Pi$ (that lives on a space $\Xi$) with respect to a vector $\uz \in  \Xi^r$ is the law of $\Pi$ conditioned to contain the points in $\Xi$ which are the co-ordinates of $\uz$. In subsequent discussions in Section \ref{supp}, we will provide a rigorous description of Palm measures. 

Let $\P_{\uz}$ denote the Palm measure of $\Pi$ with respect to $\uz$. We are interested in the mutual singularity (and continuity) of $\P_{\uza}$ and $\P_{\uzb}$ for two diffrent vectors $\uza$ and $\uzb$. According to the heurisitic description of the Palm measure (as also the rigorous definition to follow in Section \ref{supp}), a random point configuration $\xi$ sampled from the Palm measure necessarily includes the points corresponding to $\uz$, which makes the above question somewhat trivial - roughly speaking, we can decide the identity of the measure by examining whether it contains the points from $\uz_1$ or $\uz_2$. However, it is often customary to think of the Palm measure to be the law of $\l(\xi \setminus \text{ the points of }\uz \r)$. Under this identification, the question of mutual regularity becomes an interesting one, and different answers can be obtained in different natural models. 

For the Poisson process, which is the most commonly studied model of point processes, the answer to the question of mutual regularity is trivial: any two Palm measures are always mutually abosolutely continuous. This follows from the spatial independence property of the Poisson process, and is valid under mild assumptions on the intensity measure of the Poisson process (principally entailing that the intensity measure does not contain atoms).

In \cite{OsSh}, Osada and Shirai studied this question with respect to the Ginibre ensemble, which is a determinantal point process arising out of the eigenvalues of non-Hermitian random matrices. They found a very interesting behaviour: any two Palm measures of the Ginibre ensemble are mutually absolutely continuous if the lengths of the conditioning vectors are equal, and they are mutually singular otherwise. This indicates that the study of point processes with strong spatial correlation can throw up surprising answers to the question of comparing Palm measures. 

Before moving ahead, we formally state their result (Theorem 1.1 in \cite{OsSh}) below, where we denote by $\G$ the law of the Ginibre ensemle:
\begin{theorem}
\label{OsShth} 
Assume that $\ul{x} \in \C^l$ and $\ul{y} \in \C^n$. If $l = n$, then $\G_{\ul{x}}$ and $\G_{\ul{y}}$ are
mutually absolutely continuous. In addition, if $l \ne n$, then $\G_{\ul{x}}$ and $\G_{\ul{y}}$ are singular to
each other.
\end{theorem}

Results on the similar comparative behaviour of reduced Palm measures have been established by Bufetov in \cite{Bu-3} for a large class of determinantal point processes on $\R$ with integrable projection kernels (including the sine, Airy and Bessel processes), by Bufetov and Qiu for a large class of determinantal processes on $\C$ (\cite{BQ}), and for the Gamma process in 1D, a similar result goes back to the work of Olshanski (\cite{Ol}).

A key feature of these results is that they depend crucially on the determinantal structure of these models. In this work, we obtain results comparing Palm measures, similar in spirit to (and extending) Theorem \ref{OsShth}, for a wider class of point processes, particularly those not having any determinantal structure. Key examples of this include the zeroes of the standard planar Gaussian analytic function and its generalizations. Our results also exhibit more delicate dependence on the vectors $\ul{x}$ and $\ul{y}$ for absolute continuity (as contrasted with the simple dependence on dimension for the determinantal processes discussed above). 

The study of spatial conditioning in point processes with strong correlation structures has attracted a fair amount of interest in recent years. Principal examples of models studied in this regard include the Ginibre ensemble, the zeroes of the standard planar Gaussian analytic function (henceforth abbreviated as GAF), the sine kernel process on the real line, the Airy and Bessel processes, and so on. In \cite{GP}, the authors showed that in the Ginibre ensemble, the points outside a domain $D$ determine exactly the number of points in $D$. For the standard planar GAF zero process, they showed that the points outside $D$ determine the number as well as the sum of the points in $D$. Furthermore, they showed that the point configuration outside $D$ does not determine anything more about the points inside $D$. To give an idea of the precise technical sense in which these assertions hold, we will quote the relevant theorems for the Ginibre and the planar GAF zero ensembles (denoted by $\G$ and $\F$ respectively). 

In order to do so, we need to introduce some notation, which will come in handy for understanding our main results as well. A (simple) point process $\Pi$ is a random locally finite point configuration on some metric space $\Xi$ that is equipped with a regular Borel measure $\mu$. A point process can equivalently be looked upon as a random counting measure, with atoms corresponding to the points. For more details on point processes, we refer the reader to \cite{DV} and \cite{Ka}. In particular, we will be making use of the notion of the $r$-point intensity measures of a point process, for a concrete definition of which we direct the reader to \cite{HKPV} (Chapter 1, Section 1.2).

The space $\S$ of locally finite point configurations on $\Xi$ is a Polish space, and a point process $\Pi$ on $\Xi$ can be thought of as a probability measure on $\S$.  Let $\D \subset \Xi$ be a bounded open set. The decomposition $\Xi=\D \cup \D^\c$ induces a factorization $\S=\S_{\inn} \times \S_{\out}$, where $\S_{\inn}$ and $\S_{\out}$ are respectively the spaces of  finite point configurations on $\D$ and locally finite point configurations on $\D^\c$.  This immediately leads to the natural decomposition $\Upsilon=(\Upsilon_{\inn},\Upsilon_{\out})$ for any $\Upsilon \in \S$, and  consequently a decomposition of the point process $\Pi$ as $\Pi=(\Pi_{\inn},\Pi_{\out})$. 

We are now ready to state from \cite{GP} the results for the Ginibre and the planar GAF zero ensembles; in what follows, $\D$ is a bounded open set in $\C$.

\begin{theorem}
 \label{gin-1}
For the Ginibre ensemble, there is a measurable function $N:\S_{\out} \to \nat$ such that a.s. \[ \mathrm{ \ Number \ of \ points \ in \ } \G_{\inn} = N(\G_{\out})\hspace{3 pt}.\]
\end{theorem}

Since a.s. the length of $\uz$ equals $N(\G_{\out})$, we can assume that each measure $\rho(\uout,\cdot)$ is supported on $\D^{N({\uout})}$.

\begin{theorem}
\label{gin-2}
For the Ginibre ensemble, $\P$-a.s. the measure $\rho({\G_{\out}},\cdot)$ and the Lebesgue measure $\el$ on $\D^{N(\G_{\out})}$ are mutually absolutely continuous.
\end{theorem}

In the case of the GAF zero process, we prove that the points outside $\D$ determine the number as well as the centre of mass (or equivalently, the sum) of the points inside $\D$, and ``nothing more''. 

\begin{theorem}
 \label{gaf-1}
For the GAF zero ensemble, \newline
\noindent
(i)There is a measurable function $N:\S_{\out} \to \nat$ such that a.s.  \[ \mathrm{ \ Number \ of \ points \ in \ } \F_{\inn} = N(\F_{\out}).\]
(ii)There is a measurable function $S:\S_{\out} \to \C$ such that a.s.  \[ \mathrm{ \ Sum \ of \ the \ points \ in \ } \F_{\inn} = S(\F_{\out}).\]
\end{theorem}

For a possible value $\uout$ of $\F_{\out}$, define the set of admissible vectors of inside points (obtained by considering all possible orderings of such inside point configurations) \[ \Sigma_{S({\uout})} := \{ \uz \in \D^{N({\uout})} : \sum_{j=1}^{N({\uout})} \z_j = S({\uout}) \}\] where $\uz=(\z_1,\cdots,\z_{N({\uout})})$.

Since a.s. the length of $\uz$ equals $N({\uout})$, we can assume that each measure $\rho(\uout,\cdot)$ gives us the distribution of a random vector in  $\D^{N({\uout})}$ supported on $\Sigma_{S({\uout})}$.

\begin{theorem}
\label{gaf-2}
For the GAF zero ensemble, $\P$-a.s. the measure $\rho(\F_{\out},\cdot)$ and the Lebesgue measure $\els$ on $\Sigma_{S(\F_{\out})}$ are mutually absolutely continuous.
\end{theorem}

We refer to as ``rigidity'' the phenomenon in which the point configuration outside $D$ exactly determines certain statistics of the points inside $D$. By ``tolerance'', we refer to the phenomenon in which, subject to the constraints imposed by the rigidity properties, points inside $D$ can be found in ``almost any'' possible configuration.

Formally, we define rigidity and tolerance as follows. 
 \begin{definition}
 \label{defrig}
A measurable function $f_{\inn} : \S_{\inn} \rightarrow \C $ is said to be \textbf{rigid} with respect to the point process $X$ on $\S$ if there is a measurable function  $f_{\out}:\S_{\out} \rightarrow \C$ such that a.s. we have $f_{\inn}(X_{\inn})=f_{\out}(X_{\out})$.
\end{definition}
In order to give a rigorous definition of tolerance, we restrict our scope a little bit, though the present definition still captures all the known examples, and is easily amenable to generalization. 

\begin{definition}
 \label{linstat}
Let $\Pi$ be a point proces on $\Xi$ and $\ph:\Xi \to \C$ be a measurable function. Then the linear statistics $\La(\ph)[\Pi]$ is defined to be the random variable \[\La(\ph)[\Pi]:=\int_{\Xi} \ph(z) \mathrm{d}[\Pi](z),\] where $\mathrm{d}[\Pi]$ is the (random) counting measure naturally associated with $\Pi$.  
\end{definition}

\begin{definition}
 \label{tolerance}
 Let $\Pi$ be a point process on a Riemannian manifold $\Xi$ with volume measure $\mu$. Let $\D \subset \Xi$ be a bounded open set, and let $\{\La(\Phi_j)[\Pi_\inn]\}_{j=0}^t$ be rigid linear statistics, with $\Phi_0 \equiv 1$ and $\Phi_1,\cdots,\Phi_t:\D \rightarrow \C$ smooth functions. \newline
 For an integer $m \ge 0$  and $\us:=(s_1,\cdots,s_t) \subset \C^t$, consider the submanifold of $\D^m$ \[\Sigma_{m,\us}:=\{ \uz=(\zeta_1,\cdots,\zeta_m) \in \D^m : \La(\Phi_j)[\delta_{\uz}]=s_j; 1 \le j \le t  \},\] where $\delta_{\uz} \in \S$ is the point configuration corresponding to the point set $\{\zeta_i\}_{i=1}^m$.
 
 Then $\Pi$ is said to be tolerant subject to  $\{\La(\Phi_j)[\Pi_\inn]\}_{j=0}^t$ if the conditional distribution $\l(\Pi_\inn | \Pi_\out=\o \r)$ is mutually absolutely continuous with the point process of $\La(\Phi_0)[\Pi_\inn]=N(\o)$ points sampled independently from the submanifold $\Sigma_{N(\o),\us}$ (where $s_i=\La(\Phi_i)[\Pi_\inn]=S_i(\o), 1\le i \le t$) equipped with the restriction of the volume measure $\mu^{\otimes N(\o)}$.
\end{definition}

\begin{remark}
 One can generalise the above notion of tolerance by demanding constraints on more general functionals than linear statistics. For example, for a fixed positive integer $k$, one can consider a smooth function $\Psi:\D^k \to \C$ and demand that the functional $\La^k(\Psi)[\Pi_\inn]:=\int \ldots \int \Psi \mathrm{d}[\Pi_\inn]^{\otimes k}$ is rigid.  
\end{remark}

Finally, we define a \textsl{regular} collection of smooth functions:
\begin{definition}
\label{regular}
Consider a collection of smooth functions $\{\Phi_1,\cdots,\Phi_k\}$, each mapping $\Xi \to \C$. We also consider the associated function \[ \Psi_r : \Xi^r \to \C^k \] given by \[  \Psi_r(\uz):= (\La(\Phi_1)[\delta_{\uz}],\cdots,\La(\Phi_k)[\delta_{\uz}]).  \] 
We call such a collection to be regular if, for each $r \ge k$,  the Jacobian of the function $\Psi_r$ is of full rank a.e.
\end{definition}

The phenomena of rigidity and tolerance have been used to understand various questions regarding point processes, particularly those with a stochastic geometric flavour. In \cite{G12}, the rigidity of the sine kernel process was used in order to settle a natural completeness question regarding random exponential functions arising out of the sine process. More generally (Theorem 1.3 therein), a positive resolution was obtained with regard to a natural completeness question for determinantal point processes, under the assumption that the point process in question exhibits rigidity with regard to the number of points. In \cite{GKP}, the authors used rigidity and tolerance phenomena from \cite{GP} to study continuum percolation on the Ginibre and the standard planar GAF zero ensembles, in particular to establish the uniqueness of the infinite cluster in the supercritical regime. In \cite{Os}, Osada used a related quasi-Gibbs property in order to study dynamics on the Ginibre ensemble. In \cite{Bu-1} and \cite{Bu-2}, Bufetov et al. examine further interesting models of point processes from the perspective of rigidity, and obtained proofs of the rigidity of the number of points for the Airy and the Bessel processes, among others.

\section{Statement of main results}
\label{mainresults}
In this article, we explore the connections between rigidity phenomena and results of the nature of Theorem \ref{OsShth}. More specifically, we extend the results of \cite{OsSh} to point processes exhibiting rigidity and tolerance phenomena of a given nature. In heuristic terms, we show that for a point process exhibiting rigidity behaviour with respect to the statistics $\{m_i\}_{i=0}^k$ (and ``nothing more''), the Palm measures $\P_{\uza}$ and $\P_{\uzb}$ are mutually absolutely continuous if $m_i(\uza)=m_i(\uzb), 0\le i \le k$, and they are mutually singular otherwise. This fits in nicely with Theorem \ref{OsShth}, given the fact that Theorems \ref{gin-1} and \ref{gin-2} show that the Ginibre ensemble exhibits rigidity in the number of points (and ``nothing more''). 

However, it also enables us to obtain similar theorems about the mutual regularity of different Palm measures with respect to many other models, which often exhibit a much more complicated correlation behaviour than the Ginibre (already discernible in the joint density structure of the finite particle approximations). E.g., the zeroes of the standard planar GAF exhibit interactions of all orders (as opposed to pairwise interactions like in the Ginibre ensemble). 

We formally state our main theorem as follows:
\begin{theorem}
 \label{main}
 Let $\Pi$ be a point process on a Riemannian manifold $\Xi$ (without boundary) with volume measure $\mu$, and having $r$-point intensity measures $\rho_r$ mutually absolutely continuous w.r.t. $\mu^{\otimes r}$ for all $r$. Let $\Phi_0 \equiv 1$, and let $\Phi_1,\cdots,\Phi_k$ be a regular collection of smooth functions mapping $\Xi \to \C$  such that,  for any bounded open set $\D \subset \Xi$, the linear statistics $\{\La(\Phi_i)[\Pi_\inn]\}_{i=0}^k$ are rigid, and $\Pi$ is tolerant subject to $\{\La(\Phi_i)[\Pi_\inn]\}_{i=0}^k$. Then, for a.e. pair of vectors $(\uz_1,\uz_2)\in \Xi^m \times \Xi^n$, the reduced Palm measures $\P_{\uz_1}$ and $\P_{\uz_2}$ of $\Pi$ (at $\uz_1,\uz_2$ respectively) are mutually singular if \[\l(\La(\Phi_i)[\del_{\uz_1}]\r)_{i=0}^k \ne \l(\La(\Phi_i)[\del_{\uz_2}]\r)_{i=0}^k.\] Conversely, for every $r \ge k$ and a.e. $\aa \in \C^k$, $\P_{\uz_1}$ and $\P_{\uz_2}$ are mutually absolutely continuous for a.e. pair $(\uz_1,\uz_2) \in \Xi^r$ such that $\La(\Phi_0)[\del_{\uz_1}]=\La(\Phi_i)[\del_{\uz_1}]=r$ and \[\l(\La(\Phi_i)[\del_{\uz_1}]\r)_{i=1}^k = \l(\La(\Phi_i)[\del_{\uz_2}]\r)_{i=1}^k.\]  
\end{theorem}

\begin{remark}
For the singularity statement in Theorem \ref{main}, the pair $(\uz_1,\uz_2)$ are a.e. with respect to the measure $\rho_m \times \rho_n$ (equivalently, $\mu^{\otimes m} \times \mu^{\otimes n}$) on $\Xi^m \times \Xi^n$. Recall the map $\Psi$ from Definition \ref{regular}. For the absolute continuity statement, $\aa \in \C^k$ is a.e. with respect to the push-forward of $\rho_r$ (equivalently, $\mu^{\otimes r}$) under $\Psi$ and $(\uz_1,\uz_2)$ are a.e. with respect to the induced measure (from $\Xi^r$) on the submanifold \[\mathcal{M}_{\aa}=\{\uz: \l(\La(\Phi_i)[\del_{\uz}]\r)_{i=1}^k =\aa\} \subset \Xi^r.\]
\end{remark}

\begin{remark}
Theorem \ref{main} goes through verbatim (with the same proof) if each $\Phi_i$ maps into $\R$ instead of $\C$. We use $\C$ in the present article because many of our interesting examples, including the zeroes of Gaussian analytic functions, are naturally covered in that setting.
\end{remark}

One of the foremost instances where Theorem \ref{main} extends the state of the art is the case of the standard planar GAF zero process.  
In \cite{GK}, the authors introduce a family of point processes, which are zeroes of Gaussian analytic functions indexed by a parameter $\a$. These ensembles are called $\a$-GAFs, and they establish that for $\a \in (\frac{1}{k},\frac{1}{k-1}]$, the $\a$-GAF zero process exhibits rigidity at level $k$. That is, the configuration outside a domain determines the number and the first $k-1$ moments of the inside points, and ``nothing more''. Consequently, our result implies that for a.e. $\uza$ and $\uzb$, the measures $\P_{\uza}$ and $\P_{\uzb}$ are mutually absolutely continuous if the first $k$ moments of (the co-ordinates of) $\uza$ and $\uzb$ are the same, and they are mutually singular otherwise. This shows, in particular, that the mutual regularity properties of the different Palm measures of a point process can depend on the conditioning vector in an arbitrarily complicated manner (the complexity of the dependence being measured by the number of  statistics that need to be matched in order to ensure absolute continuity).

In \cite{HoOs} and \cite{OsTa}, a \textsl{quasi Gibbs} property is established for sine, Airy ($\beta=1,2,4$)  and Bessel ($\beta=2$) point processes. Rigidity of the number of points for these processes was established in \cite{Bu-1}, and tolerance of these point processes (subject to the number of points) can be deduced from this quasi Gibbs property. Consequently, we can invoke our Theorem  \ref{main} to obtain a new proof of the analogue of Theorem \ref{OsShth} for these processes.

We formally state these results as follows (we denote by $\uz(i)$  the $i$-th co-ordinate of the vector $\uz$):
\begin{theorem}
 \label{instances}
 Let $\Pi$ be a point process on a Riemannian manifold $\Xi$ with volume measure $\mu$, and let $\uz_1 \in \Xi^m,\uz_2 \in \Xi^n$. Then the following statements are true about the reduced Palm measures $\P_{\uz_1}$ and $\P_{\uz_2}$ :
 \begin{itemize}
  \item When $\Pi$ the standard planar GAF zero process on $\C$,  
        \begin{itemize}
               \item  For Lebesgue-a.e. $s \in \C$ and a.e. $\uz_1,\uz_2$ such that $|\uz_1|=|\uz_2|$ and          $\sum_{i=1}^{|\uz_1|} \uz_1(i) = \sum_{i=1}^{|\uz_2|} \uz_2(i) =s$,  $\P_{\uz_1}$ and $\P_{\uz_2}$ are mutually absolutely continuous.
               \item For a.e. $\uz_1,\uz_2$ such that $|\uz_1| \ne |\uz_2|$ or $\sum_{i=1}^{|\uz_1|} \uz_1(i) \ne \sum_{i=1}^{|\uz_2|} \uz_2(i) $, $\P_{\uz_1}$ and $\P_{\uz_2}$ are mutually singular.
        \end{itemize}
  \item When $\Pi$ the $\a$-GAF zero process on $\C$,  
        \begin{itemize}
               \item  For Lebesgue-a.e. $\m \in \C^{\lfloor \frac{1}{\a} \rfloor}$ and a.e. $\uz_1,\uz_2$ such that $|\uz_1|=|\uz_2|$ and          $\sum_{i=1}^{|\uz_1|} \uz_1(i)^j = \sum_{i=1}^{|\uz_2|} \uz_2(i)^j =\m(j)$ for all $1\le j \le \lfloor \frac{1}{\a} \rfloor$,  $\P_{\uz_1}$ and $\P_{\uz_2}$ are mutually absolutely continuous.
               \item For a.e. $\uz_1,\uz_2$ such that $|\uz_1| \ne |\uz_2|$ or $\sum_{i=1}^{|\uz_1|} \uz_1(i)^j \ne \sum_{i=1}^{|\uz_2|} \uz_2(i)^j $ for some $1\le j \le \lfloor \frac{1}{\a} \rfloor$, $\P_{\uz_1}$ and $\P_{\uz_2}$ are mutually singular.
        \end{itemize}
 \item For $\Pi$ the i.i.d. perturbation of $\Z^2$ (resp., $\Z$) with random variables having a non-vanishing density on $\R^2$ (resp. $\R$) with a finite second (resp., first) moment, we have, for Lebesgue-a.e. $(\uza,\uzb) \in \C^m \times \C^n$ (resp., $\R^m \times \R^n$),   
  \begin {itemize}
 \item   $\P_{\uz_1}$ and $\P_{\uz_2}$ are mutually singular  if $m \ne n$
 \item $\P_{\uz_1}$ and $\P_{\uz_2}$ are mutually absolutely continuous  if $m = n$
\end{itemize}
 \end{itemize}
\end{theorem}
\begin{remark}
For i.i.d. perturbations of $\Z^d$ by $d$-dimensional Gaussians having small enough variance, a similar conclusion as the 1 or 2D lattice perturbations above holds.
\end{remark}

This theorem  follows from our main Theorem \ref{main}, coupled with the results  on the rigidity and tolerance properties of these ensembles, as in \cite{GP} (Theorem 1.1 - Theorem 1.4) and \cite{GK}  (Theorem 2.1). Rigidity of the number of points for i.i.d. lattice perturbations satisfying the conditions in the statement of Theorem \ref{instances} has been established in \cite{PS}; the tolerance (in our terminology) is a simple consequence of the fact that the perturbations are independent and have a positive density a.e. with respect to the  Lebesgue measure. The remark about Gaussian perturbations also follows from rigidity established in \cite{PS} and a tolerance statement that follows from independence considerations.

Here we illustrate the details in the case of the standard planar GAF zero process; the details in the other cases are on similar lines. Our goal is to verify that  the standard planar GAF zero process satisfies the conditions of Theorem \ref{main}. For this, we will make use of the rigidity-tolerance behaviour of this point process, which was established in \cite{GP}; for convenience the relevant results have been quoted here as Theorems \ref{gaf-1} and \ref{gaf-2}. In terms of the conditions laid out in Theorem \ref{main}, clearly $\Xi = \C$ in this case, with $\mu$ the Lebesgue measure on $\C$. It is well known (also easy to see from the definition of the GAF) that the $r$-point intensity measures of the GAF zeros have densities with respect to the Lebesgue measure on $\C^r$ (for a specific reference, see \cite{HKPV}). We put $k=1$ and $\Phi_1(z)=z$. Clearly, the function $\Phi_1$ is a \textsl{regular} as per Definition \ref{regular}, as can be seen from the fact that the Jacobian of the map \[\Psi_r : (z_1,\cdots,z_r) \mapsto \sum_{i=1}^r z_i\] is $[1,1,\cdots,1]$, which is always of full rank.  Theorem \ref{gaf-1} and \ref{gaf-2} are equivalent to the statement that for any bounded open set $D \subset \C$, the statistics $(\La(\Phi_i)[\Pi_\inn])_{i=0}^1=(N(\Pi_\inn),S(\Pi_\inn))$ are rigid, and the GAF zero process is tolerant subject to \newline $(N(\Pi_\inn),S(\Pi_\inn))$ (recall the Definitions \ref{defrig} and \ref{tolerance} of rigidity and tolerance respectively). This verifies the conditions of Theorem \ref{main} for the GAF zero process.

We now interpret the conclusions of Theorem \ref{main} for the GAF zero process. For a vector $\uz \in \C^r$, denote by $|\uz|$ the dimension and by $S(\uz)$ the sum of the co-ordinates of $\uz$. 
Then Theorem \ref{main} implies that for a.e.-pair $(\uz_1,\uz_2) \in \C^m \times \C^n$ such that $(|\uz_1|,S(\uz_1)) \ne (|\uz_2|,S(\uz_2))$, the reduced Palm measures $\P_\uza$ and $\P_\uzb$ are mutually singular. In particular, this implies that if $m \ne n$, then $\P_\uza$ and $\P_{\uzb}$ are mutually singular for Lebesgue a.e.-$(\uza,\uzb) \in \C^m \times \C^n$.This brings us to the situation $m = n$. In this scenario, there are two possibilities: $S(\uza)\ne S(\uzb)$ and $S(\uza) = S(\uzb)$. Regarding the former possibility, for Lebesgue a.e.$(\uza,\uzb) \in \C^m \times \C^m$, Theorem \ref{main} states that $\P_{\uza}$ and $\P_{\uzb}$ are mutually singular.  A particularly interesting case of this is when $m=n=1$, which we state as a separate corollary:
\begin{corollary}
\label{gef1}
For the standard planar GAF zero process, the reduced Palm measures $\P_z$ and $\P_w$ are mutually singular for Lebesgue a.e. pair $(z,w) \in \C \times \C$.
\end{corollary}
This contrasts markedly with the analogous comparison of Palm measures in most determinantal processes, including the Ginibre process (Theorem \ref{OsShth}).

This leaves us with the final case: $m=n$ and $S(\uza)=S(\uzb)$. Denoting $S(\uza)=S(\uzb)=s \in \C$, we consider the manifold \[\M_s :=\{ \uz \in  \C^m : S(\uz)=s         \}. \] $\M_s$ carries a natural Lebesgue measure, induced from the Lebesgue measure on $\C^m$, denote this measure by $l_s$. Then Theorem \ref{main} says that for Lebesgue a.e.-$s$, we have that for $l_s$-a.e. pair ($\uza,\uzb$) $\in \M_s$, the reduced Palm measures $\P_\uza$ and $\P_\uzb$ are mutually absolutely continuous.



\section{Proof of Theorem \ref{main}}

In this section, we prove our main Theorem \ref{main}. We will first analyze the support properties of Palm measures, and then connect it with rigidity properties of point processes, in two subsections.

\subsection{Palm measures and their support}
\label{supp}
Let, as before, $\S$ denote the Polish space of all locally finite point configurations on $\Xi$, and $\mathcal{B}(\Xi)$ denote the Borel sigma field on $\S$. We begin with the 1-point Campbell measure and a rigorous definition of the 1-point Palm measure of point processes, we refer the reader to \cite{Ka} Chapter 10 for a more detailed treatment than we present here. The 1-point Campbell measure $\mu^{(1)}$ of a point process $\Pi$ (whose law we denote by $\P$) is the measure defined on $\Xi \times \S$ given by \[\int f(s,\xi) \mathrm{d}\mu^{(1)}(s,\xi):= \int \l( \int f(s,\xi) \mathrm{d}[\xi](s) \r) \mathrm{d}\P(\xi).\] In the above equation, the measure $[\xi]$ is the counting measure that naturally corresponds to $\xi \in \S$. The 1-point Palm measures  $\{\ol{\P}_s: s \in \S\}$ (that include the points in the conditioning vector) are defined by a decomposition of the measure $\mu^{(1)}$ into a regular conditional measure with respect to the first co-ordinate:
\[\mathrm{d}\mu^{(1)}(s,\xi)=\mathrm{d}\rho_1(s) \times \mathrm{d}\ol{\P}_s(\xi). \] In other words, we define the 1-point Palm measure by the integral formulation \[ \int f(s,\xi) \mathrm{d}\mu^{(1)}(s,\xi) = \int \l( \int f(s,\xi) \mathrm{d}\ol{\P}_s(\xi) \r) \mathrm{d}\rho_1(s). \] Here $\rho_1$ is the 1-point intensity measure of $\Pi$.  

For any simple locally finite point configuration $\xi \in \S$ and an integer $r \ge 1$, let us denote by $[\xi]^{\wedge r}$ the counting measure on all possible ordered $r$-tuples of distinct points of $\xi$. Then, for any integer $r \ge 1$, the $r$-point Campbell measure $\mu^{(r)}$ can be defined as a measure on $\Xi^r \times \S$ given by \[\int f(\us,\xi) \mathrm{d}\mu^{(r)}(s,\xi):= \int \l( \int f(\us,\xi) \mathrm{d}[\xi]^{\wedge r}(\us) \r)\mathrm{d}\P(\xi).\] 
Consequently, one can define $\ol{\P}_{\us}$, the $r$-point Palm measure at $\us$ (that includes the points in $\us$) by \[\mathrm{d}\mu^{(r)}(\us,\xi)=\mathrm{d}\rho_r(\us) \times \mathrm{d}\ol{\P}_{\us}(\xi), \] or equivalently, \[\int f(\us,\xi) \mathrm{d}\mu^{(r)}(\us,\xi) = \int \l( \int f(\us,\xi) \mathrm{d}\ol{\P}_{\us}(\xi) \r) \mathrm{d}\rho_r(\us), \] where $\rho_r$ is the $r$-point intensity measure of $\Pi$. 

Since $\pi(\uz) \subset \xi$ for each $\xi \in \Supp(\ol{\P}_{\uz})$, therefore, we can equivalently consider the law of $\xi \setminus \uz$. We call this measure the \textsl{reduced} Palm measure of $\Pi$ at $\zz$. We will denote this measure by $\P_{\uz}$. 

Let $Q$ be a countable dense subset of $\Xi$. We will call a subset $G$ of $\Xi$ to be \textsl{good} if  $G$ is the union of finitely many disjoint open balls with distinct centres in $Q$ and identical rational radius. We will say that a nested sequence $\{G_n\}$ of good subsets of $\Xi$ (having a fixed number $m$ of constituent balls) \textsl{converge} to  $\underline{p} \in \Xi^{m}$ if $G_{n+1} \subset G_n$ and the centres of the constituent balls of $G_n$ converge to  $\underline{p}$ (in some ordering of the co-ordinates). In such a situation, we will say that $\underline{p} \in \Xi^m$ is a \textsl{limit} of $\{G_n\}$. Finally, we will say that $\uz \in \Xi^r, r \le m$ (with distinct co-ordinates) \textsl{belongs} to the limit $\underline{p}$ of such a  sequence ${G_n}$ of good sets (equivalently, we say that $\underline{p}$ \textsl{contains} $\uz$ ) if the co-ordinates of $\uz$ are a subset of those of $\underline{p}$.

Recall, from Section \ref{mainresults}, the notation that $\m(\uz)$ denotes the vector $(m_1(\uz),\cdots,m_k(\uz))$, where $m_i(\uz)=\Lambda(\Phi_i)[\del_{\uz}]$.  For any bounded open set $D \subset \Xi$, the number of points of our process $\Pi$ that lie in $D$ will be denoted by $N(D)$. This quantity is a measurable function of the point configuration $\theta$ in $\D^\complement$  because of the rigidity of the number of points (which corresponds to the functional $\Psi_0 \equiv 1$), and we will denote this function by $m_0(\theta;D)$.
Similarly, $\m(\Pi_\inn)$ (where $\Pi_\inn$ denotes the points of $\Pi$ inside $D$ in uniform random order) is a measurable function of $\theta$, and we denote this by $\mm(\theta;D)$.

Let $r \le p$ be positive integers. For a good set $G \subset \Xi$ having $p$ constituent balls, we define the event $\A(G,r) \in \mathcal{B}(\Xi)$, which entails that $\xi \in \A(G,r)$ if 
\begin{itemize}
	\item $\xi$ is supported on $G^\complement$
	\item $m_0(\xi;G)=r$.
\end{itemize}

For every $\uz \in G^r$, we define the event $\A(G,\uz) \in \mathcal{B}(\Xi)$, which entails that $\xi \in \A(G,\uz)$ if 
\begin{itemize}
	\item $\xi$ is supported on $G^\complement$
    \item  $m_0(\xi;G)=|\uz|=r$ 	
	\item $\mm(\xi;G)=\m(\uz)$.
\end{itemize}

Finally, we say that a sequence of sets $\{B_n\}$ exhausts the support of a measure $\gamma$ if $\gamma(B_n^\c) \to 0$ as $n \to \infty$.

With these definitions, we are ready to state the following technical result:
\begin{lemma}
	\label{second}
\begin{itemize}	
 For $\rho_r$-a.e. $\eta$, it holds that, for any nested sequence $\{G_n\}$ of good sets (having $p \ge r$ constituent balls for each $n$) with a limit that contains $\eta$ :
\item (i) The events $\A(G_n,r)$ exhaust the support of $\P_{\eta}$. 
\item (ii) The events $\A(G_n,\eta)$ exhaust the support of $\P_{\eta}$.

\end{itemize}
\end{lemma}

\begin{remark}
Since we assume that the sequence of good sets $G_n$ is nested, therefore $\eta$ is contained in the limit of $\{G_n\}$ implies that $\eta \subset G_n$ for each $n$. 
\end{remark}

\begin{proof}
	Observe that we trivially have the inclusion of events $\A(G_n,\uz) \subset \A(G_n,r)$. Therefore, it suffices to establish part (ii) of the Lemma, from which part (i) will follow.
	
	We proceed as follows. 
	First of all, for a good set $G$ (with $p$ constituent balls) and $\uz \in G^r$, consider the event $\mathcal{F}(G,\uz)$ such that a point configuration $\xi \in \mathcal{F}(G,\uz)$ entails that 
	\begin{itemize}
		\item $\xi$ is supported on $G^\complement$
		\item $(m_0(\xi;G),\mm(\xi;G)) \ne (m_0(\uz),\m(\uz))$.
	\end{itemize} 
	We assert that for $\rho_r$-a.e. $\uz \in G^r$, we have  $\P_{\uz}[\mathcal{F}(G,\uz)]=0$. To this end, 
	we observe that
	\begin{align*}
	&\int_{G^r}       \P_{\uz}[\mathcal{F}(G,\uz)]              \rho_r(\uz) \d V(\uz)  \\
	=& \P[\uz \cup \xi  \text{ is a realisation of } \Pi \text{ for some } \uz  \in G^r \text{ and some }  \xi \in \mathcal{F}(G,\uz)]   \\
	=& 0, 
	\end{align*}
	where, in the last step we have used the fact that  \[ \P[\uz \cup \xi  \text{ is a realisation of } \Pi \text{ for some } \uz  \in G^r \text{ and some }  \xi \in \mathcal{F}(G,\uz)] =0 \] because of the rigidity properties of $\Pi$ with respect to the set $G$. More precisely, since $\{\uz \cup \xi \} \cap G =\uz $ and $\{\uz \cup \xi \} \cap G^\c =\xi $ (as point sets), therefore by the rigidity of $\Pi$ we have $(m_0(\xi;G),\mm(\xi;G)) $ must equal $(m_0(\uz),\m(\uz))$ (for $\P$-a.e. realisation $\omega=\uz \cup \xi$ of the point process such that $\omega \cap G = \uz$ and $\omega \cap G^\complement =\xi$). This proves the assertion. 
	
	Since there are only countably many good sets, we can deduce from the above that for $\rho_r$-a.e. $\uz$, we have  $\P_{\uz}[\mathcal{F}(G,\uz)]=0$ for any good set $G$ such that  $\uz \in G^r$.
	
	Now let us consider a $\uz$ satisfying the above assertion, and a nested sequence of good sets $\{G_n\}$ (with $p$ constituent balls each) having a limit that contains $\uz$ (and, consequently, $\uz \subset G_n$ for each $n$). Consider the event $\A(G_n,\uz)^\complement$, under the reduced Palm measure $\P_{\uz}$. This event can occur only in two ways (respectively corresponding to the  defining conditions of the event $\A(G_n,\uz)$): 
	\begin{itemize}
        \item There is at least one point of the Palm process $\P_{\uz}$ inside $G_n$. 		
		\item $\mathcal{F}(G_n,\uz)$ occurs
	\end{itemize}
	By choice of $\uz$, we already have $\P_{\uz}[\mathcal{F}(G_n,\uz)]=0$. Thus, recalling  that $N(U)$ denotes the number of points of a configuration that lie in the set $U$, we have
	\[\P_{\uz}[\A(G_n,\uz)^\complement] \le \P_{\uz}[N(G_n) \ge 1] \le \E_{\uz}[N(G_n))] \downarrow 0 \] as $n \to \infty$, by the Dominated Convergence Theorem.

\end{proof}

\subsection{Palm measures and rigidity phenomena}

\subsubsection{Singularity}
\label{singularity}

For $\uz \in \Xi^p, p \ge 1$, recall the notation  $m_i(\uz)=\La(\Phi_i)[\del_{\uz}], 0 \le i \le k$, and $\m(\uz)=(m_1(\uz),\cdots,m_k(\uz))$. Also recall that $m_0(\uz)=|\uz|$.

Consider $\uza \in \C^r$ and $\uzb \in \C^s$ such that $(m_0(\uza),\m(\uza)) \ne (m_0(\uzb),\m(\uzb))$, both satisfying the conclusions of Lemma \ref{second} part (ii) (this happens a.e.-$\rho_r \times \rho_s$). We also assume that $\uza$ and $\uzb$ have distinct co-ordinates (both within and between themselves), since this also happens a.e.-$\rho_r \times \rho_s$. Let $\{G_n\}$ be a nested sequence of good sets, each having $r+s$ constituent balls (and each ball containing exactly one co-ordinate of either $\uza$ or $\uzb$), such that $\{G_n\}$ has the limit $(\uza,\uzb)$ (in the sense of the definitions in the previous subsection).

By Lemma \ref{second} part (ii), the support of $\P_\uza$ is exhausted by $\A(G_n,\uza)$ and the support of $\P_\uzb$ is exhausted by $\A(G_n,\uzb)$. But since $(m_0(\uza),\m(\uza)) \ne (m_0(\uzb),\m(\uzb))$, therefore $\A(G_n,\uza) \cap \A(G_n,\uzb) = \phi$. In other words, $\A(G_n,\uza) \subset \A(G_n,\uzb)^\c$ and $\A(G_n,\uzb) \subset \A(G_n,\uza)^\c$. 

We make the following claim :  for two probability measures $\mu_1$ and $\mu_2$ on the same space, suppose there is a sequence of events $B_n$ such that for any $\eps >0$, $\mu_1(B_n)>1-\eps$ for large enough $n$, and $\mu_2(B_n)<\eps$ for large enough $n$. Then $\mu_1$ and $\mu_2$ are mutually singular.


Before proving this claim, we note that this suffices to complete the proof of singularity. To see this, set $\mu_1=\P_\uza$, $\mu_2=\P_\uzb$ and $B_n=\A(G_n,\uza)$. For any $ \eps>0$ we note that $\P_\uza(\A(G_n,\uza))>1-\eps$ for all large enough $n$ because these sets exhaust the support of $\P_\uza$. But $\A(G_n,\uza) \subset \A(G_n,\uzb)^\c$, so $\P_\uzb (\A(G_n,\uza))<\eps$ for all large enough $n$, because $\A(G_n,\uzb)$-s exhaust the support of $\P_\uzb$. Then, from the above claim, it follows that $\P_\uza$ and $\P_\uzb$ are mutually singular.

It remains to prove the claim. Let $\mu_1,\mu_2,\{B_n\}$ be as in the claim.  Passing to a sub-sequence if necessary, we may assume that $\sum_n \mu_2(B_n) <\infty$. Consider the event \[ B:= \varlimsup B_n := \cap_{N=1}^{\infty} \cup_{n\ge N} B_n . \] Let $C_N$ denote $\cup_{n\ge N} B_n$. For $\eps>0$ and $N$ large enough, $\mu_1(C_N) \ge \mu_1(B_{N}) > 1-\eps$. But  the $C_N$-s are decreasing in $N$, and hence $\mu_1(B)= \lim_{N \to \infty} \mu_1(C_N)=1$. On the other hand, \[\mu_2(C_N) \le \sum_{n \ge N} \mu_2(B_n).  \] Since $\sum_n \mu_2(B_n) <\infty$, therefore the right hand side can be made arbitrarily small by choosing $N$ large enough. Hence we have $\mu_2(B)=\lim_{N \to \infty} \mu_2(C_N) =0$. 

The upshot of this is that $\mu_1(B)=1$, whereas $\mu_2(B)=0$. Since $\mu_1,\mu_2$ are probability measures, this implies that  $\mu_1(B^\c)=0$ and $\mu_2(\B^\c)=1$. This completes the proof that the measures $\mu_1$ and $\mu_2$ are mutually singular.

This completes the proof that of mutual singularity of $\P_\uza$ and $\P_{\uzb}$.

\subsubsection{Absolute Continuity}
\label{absconti}
For $\aa \in \C^r$, denote by $\M_{\aa}$ the set $\uz \in \Xi^r$ such that $\m(\uz)=\aa$. Consider the $r$-point intensity measure $ \rho_r \d V$ on $\Xi^r$ (where $\d V$ is the canonical volume measure on $\Xi^r$, and $\rho_r$ is the $r$-point intensity function).  Consider the map $\Psi_r: \Xi^r \mapsto \C^k$ given by $\uz \to \m(\uz) $, which is of full rank . This implies that we can decompose $\rho_r \d V$  (going to local co-ordinates if necessary) as $  \mu(\aa)\d \aa \times \nu(\aa,\uz) \d l_{\aa}(\uz)$, where $\d \aa$ is Lebesgue measure on $\C^k$, $\mu(\aa) \d \aa$ is the push forward of $\rho_r \d V$ to $\C^k$ under $\Psi_r$,  $\d l_{\aa}$ is the induced measure on $\M_{\aa}$ from $\d V$, and $\nu(\aa,\uz)$ is a  density (see, e.g., \textsl{smooth co-area formula}, \cite{Ch}, Chap. III). Roughly speaking, this corresponds to a foliation of $\Xi^r$ by the level sets of $\m$. Since $\Xi$ is covered by a countable union of such charts (i.e., local neighbourhoods), therefore it suffices to work on each such chart. 

Since $\rho_r > 0$ a.e. with respect to the canonical volume measure of $\Xi^r$, therefore for $\mu(\aa)\d \aa$-a.e. $\aa$, we have $\nu(\aa,\uz)>0$ for $\d l_{\aa}$-a.e. $\uz$. Consequently,  we deduce the following

\begin{claim}
	\label{first}
For $\mu(\aa)\d \aa$-a.e. $\aa$, it is true that for $\d l_{\aa} \times \d l_{\aa}$-a.e. pair $(\uz_1,\uz_2)$ (so that $\m(\uz_1)=\m(\uz_2)=\aa)$, we have  $\nu(\aa,\uz_1)$,$\nu(\aa,\uz_2)>0$.
\end{claim}

In light of Lemma \ref{second}, we deduce that for  $\mu(\aa)\d \aa$-a.e. $\aa$, we have that $\d l_{\aa}$-a.e. $\uz$ satisfies the support properties as in the conclusion of Lemma \ref{second}.

Let $\uz_1$ and $\uz_2$ be two configurations such that $\m(\uz_1)=\m(\uz_2)$ and they satisfy the properties laid out in Claim \ref{first} and Lemma \ref{second}. 
Let $\{D(\eps)\}_{\eps \downarrow 0} \subset \Xi$ be a nested sequence of good sets (heuristically speaking, they approximate the $\eps$ neighbourhood of the co-ordinates of $\uz_1$ and $\uz_2$) such that their limit is the $2r$ dimensional configuration $(\uz_1,\uz_2)$.   

For brevity, in what follows, we will denote by  $\mathcal{A}_\eps$   the event $\mathcal{A}(D(\eps),r)$ (as defined in Section \ref{supp}). Recall that, by definition,  each point configuration $\xi \in \mathcal{A}_\eps$ satisfies the following conditions:
\begin{itemize}
	\item $\xi$ is supported on $D(\eps)^\complement$ 
	\item $m_0(\xi;D(\eps))=r$.
\end{itemize}

We consider the joint law  $(\Pi_\inn,\Pi_\out)$ of the points in $D(\eps)$ and  in $D(\eps)^\complement$ respectively. We denote by $E_r$ the event that there are $r$ points in $D(\eps)$; equivalently $|\Pi_\inn|=r$. By the rigidity of the number of points in $D(\eps)$, the event $E_r$ is measurable with respect to $\Pi_\out$.  We consider the law of $(\Pi_\inn,\Pi_\out)$ conditioned on the event $E_r$. The rigidity properties of $\Pi$ imply that, on the event $E_r$, the random variable $\m(\Pi_\inn)$  is measurable with respect to  $\Pi_\out$. This implies that there is a regular conditional probability $Q_\eps(\m(\Pi_\inn), d \xi)$ that pertains to the random variable $\Pi_\out$ given  $\m(\Pi_\inn)$ and given that $E_r$ occurs.

 For a given vector $\aa \in \C^k$, let $\Sigma_{\aa}$ denote the subset of $D(\eps)^r$ such that $\m(\uz)=\aa$ for all $\uz \in \Sigma_{\aa}$; clearly $\Sigma_{\aa} = D(\eps)^r \cap \M_{\aa}$.
From the rigidity properties of our point process, we know that the random variable $\m(\Pi_\inn)$ is a measurable function of $\Pi_\out$. Moreover, we also know that on the event $E_r$ (measurable with respect to $\Pi_\out$),  the conditional law of $\Pi_\inn$ given $\Pi_\out$, denoted  $ \d\P(\uz|\Pi_\out=\xi)$, has a density  $f(\uz,\xi)$  with respect to the measure $\d l_{\m(\uz)}$ on $\Sigma_{\m(\uz)}$ (note here that $\Sigma_{\m(\uz)}$ is determined by $\xi$ because $\m(\Pi_\inn)$ is measurable with respect to $\Pi_\out$).

We now consider the reduced Palm measure with respect to configurations in $D(\eps)^r$ (on the event $\mathcal{A}_\eps$).  In other words, we consider the measures $\P_{\uz}(\cdot \cap \mathcal{A}_\eps)$, where $\uz \in D(\eps)^r$. To introduce our candidate  for  $\P_{\uz}(\cdot \cap \mathcal{A}_\eps)$, we need to first introduce another quantity. 

We can consider the measure $\kappa$ on $\D^r$ which is the marginal distribution of the points of $\Pi$ inside $\D$ on the event $E_r$ (with the points being taken in uniform random order). In other words, this is the measure $\P[E_r] \d(\Pi_\inn | E_r)$, where $(\Pi_\inn|E_r)$ is the law of $\Pi_\inn$ given $E_r$ occurs.  Clearly, this measure is absolutely continuous with respect to the measure $\rho_r \d V$. Consequently, the push forward $(\m)_* \kappa$ is absolutely continuous  with respect to the push forward $(\m)_* [\rho_r \d V]$. The latter measure, as we may recall is $\mu(\aa) \d \aa$, which means that there exists a density $\varrho$  such that 
\begin{equation}
\label{marginal1}
[ (\m)_* \d \kappa] (\aa) = \varrho(\aa) \mu(\aa) \d \aa.  
 \end{equation}
Also, it follows from the above discussion that
\begin{equation}
\label{marginal2}
\P[E_r] [ (\m)_* \d  (\Pi_\inn |E_r)] (\aa) = [ (\m)_* \d \kappa] (\aa).
\end{equation}

Our candidate for  $\P_{\uz}(\cdot \cap \mathcal{A}_\eps)$ (for $\uz \in D(\eps)^r$) is the following: \begin{equation} \label{kap} \beta_{\uz}( \d \xi) = \varrho(\m(\uz))  f(\uz,\xi) Q(\m(\uz),\d \xi)/ \nu(\m(\uz),\uz), \end{equation} and for  $\uz$ such that $\nu(\m(\uz),\uz)=0$ we simply define $\beta_{\uz}( d \xi)$ to be 0. In other words, we claim that $\P_{\uz}(\cdot \cap \mathcal{A}_\eps) = \beta_{\uz}(\cdot \cap \mathcal{A}_\eps)$ for a.e. $\uz \in D(\eps)^r$. 

Let us check that this indeed true for a.e. $\uz \in D^r$.
Consider an event $ U \times (A \cap \mathcal{A}_\eps)$ such that $U \subset D(\eps)^r$ and $A$ is a measurable set in $\mathcal{B}(\Xi)$.  We set $A'=A \cap \mathcal{A}_\eps$, and denote by $M_\eps \subset \C^k$ the image of $D(\eps)^r$ under the map $\underline{m}$.

We have,
\begin{align*} 
&\int_U   \rho_r(\uz) \l( \int_{A'} \beta_{\uz}(\d \xi) \r) \d V(\uz) \\ 
=& \int_{M_\eps} \mu(\aa) \l( {\int_{\Sigma_{\aa} \cap U}} \nu(\aa,\uz) \l(\int_{A'} {\beta_{\uz}} (d\xi) \r)   \d l_{\aa} (\uz) \r)               \d\aa \\
=& \int_{M_\eps} \mu(\aa) \l( {\int_{\Sigma_{\aa} \cap U} }  \l(\int_{A'}  \nu(\aa,\uz) {\beta_{\uz}} (d\xi) \r)   
\d l_{\aa} (\uz) \r)               \d\aa \\
=& \int_{M_\eps} \mu(\aa) \l( {\int_{\Sigma_{\aa} \cap U} }  \l(\int_{A'}  \nu(\aa,\uz) \varrho(\aa) f(\uz,\xi) Q_\eps(\aa,\d \xi)/ \nu(\aa,\uz) \r)    \d l_{\aa} (\uz) \r)  \d\aa  \\
=& \int_{M_\eps} \varrho(\aa) \mu(\aa) \l( {\int_{\Sigma_{\aa} \cap U}} \l( \int_{A'}   f(\uz,\xi) Q_\eps(\aa,\d \xi)   \r) \d l_{\aa} (\uz) \r)  \d\aa  \\
=& \int_{M_\eps} \varrho (\aa) \mu(\aa) \l(   \int_{A'} \l( \int_{{\Sigma_{\aa}} \cap U} f(\uz,\xi)  \d l_{\aa} (\uz)   \r) Q_\eps(\aa,\d \xi) \r)  \d\aa \text{ \{by Fubini's Theorem \} }  \\
=&  \int_{M_\eps} \varrho (\aa)\mu(\aa) \l(   \int_{A'}  \P[\Pi_\inn \in U | \Pi_\out=\xi ]  Q_\eps(\aa,\d \xi) \r)  \d\aa   \text{ \{by definition of } f \} 
\end{align*}
\begin{align*} 
=&  \int_{M_\eps}  \varrho(\aa) \mu(\aa)  \P[\Pi_\inn \in U, \Pi_\out \in A' | \m(\Pi_\inn) = \aa  \cap E_r]   \d\aa   \text{ \{by definition of } Q_\eps \}  \\
=&  \int_{M_\eps}    \P[\Pi_\inn \in U, \Pi_\out \in A' | \m(\Pi_\inn) = \aa \cap E_r ]   [(\m)_*\d\kappa](\aa)   \text{ \{by \eqref{marginal1}\} } \\
=&  \int_{M_\eps}    \P[\Pi_\inn \in U, \Pi_\out \in A' | \m(\Pi_\inn) = \aa \cap E_r ] \P[E_r]   [(\m)_*\d(\Pi_\inn | E_r)](\aa)   \text{ \{by \eqref{marginal2} \} }  \\
=&  \P[E_r] \P[\Pi_\inn \in U, \Pi_\out \in A' | E_r] \text{ \{by definition of } \d(\Pi_\inn | E_r)\}  \\
=&  \P[\Pi_\inn \in U, \Pi_\out \in A' ] \text{ \{ since } \P[\Pi_\inn \in U, \Pi_\out \in A' |E_r^\c]=0 \}  \\ 
=&  \int_U   \rho_r(\uz) \l( \int_{A'} \P_{\uz}(\d \xi) \r) \d V(\uz) \text{ \{by definition of } \P_{\uz} \}. 
\end{align*}

This shows that,  for a.e. $\uz \in D(\eps)^r$,  $\P_{\uz}(\cdot \cap \A_\eps)= \beta_{\uz}(\cdot \cap \A_\eps)$.  The definition \eqref{kap} of $\beta_{\uz}$ implies that, for $\d l_{\aa}$-a.e.  $\uz,\uz' \in \Sigma_{\aa}$, we have $\beta_{\uz} \equiv \beta_{\uz'} $. Since $\uz_1$ and $\uz_2$ belong to $\Sigma_{\aa}$ for the same $\aa$, we deduce that 
\begin{equation} \label{equive} \P_{\uz_1}(\cdot \cap \A_\eps) \equiv \P_{\uz_2}(\cdot \cap \mathcal{A}_\eps). \end{equation}

For any event $B \in \mathcal{B}(\Xi)$ such that $\P_{\uz_1}(B)=0$, we have $\P_{\uz_1}(B \cap \mathcal{A}_\eps) =0$. This implies that $\P_{\uz_2}(B \cap \mathcal{A}_\eps) =0$, by the mutual absolute continuity of the  measures in \eqref{equive}. But, as $\eps \to 0$, we have $\mathcal{A}_\eps$ exhausts the support of both $\P_{\uz_1}$ and $\P_{\uz_2}$ (because $\uz_1$ and $\uz_2$ were both chosen to satisfy Lemma \ref{second}). Letting $\eps \downarrow 0$, we deduce that $\P_{\uz_2}(B)=0$. This shows that $\P_{\uz_1} \equiv \P_{\uz_2}$.

\section{Extensions}

Theorem \ref{main} can be extended in several directions. One immediate direction is the case when $\Xi$ is a countable discrete set, e.g.  $\Z, \Z^d$ or a subset thereof. In many ways, this situation is technically simpler than the continuum setting in which Theorem \ref{main} is stated and proved.  The volume measure on $\Xi$ (and its subsets) will be naturally replaced by the counting measure on those sets, and there would be no regularity assumptions on the functions $\Phi$. The proof would be the same as the proof of Theorem \ref{main}.

Another pertinent question to ask is about the situation when $|\uz_1| = |\uz_2|,  \m(\uz_1) \ne \m(\uz_2)$, but some subsets of coordinates of these two vectors match. This set will be of zero $\rho_r$ measure, and hence is not covered by Theorem \ref{main} as is. However, under a mild regularity assumption on $(\Phi_1,\cdots,\Phi_k)$, we can deal with this sceanrio as well. Recall that the assumption in Theorem \ref{main} is that the functions  $(\Phi_1,\cdots,\Phi_k)$ constitute a \textsl{regular} collection of functions (in particular, recall Definition \ref{regular}). To address the finer question, we make the additional assumption that $(\Phi_i)_{i \in S}$ is a regular collection of functions for each subset $S \subset \{1,\cdots,k\}$. Consider, for $\uz \in \Xi^r$, the vector $s(\uz)$ given by the co-ordinates of $\m(\uz)$ whose indices are in $S$. Also consider, for ${\aa} \in \C^{|S|}$, the sub-manifold $\M_{\aa}^{[S]}$ formed by $\uz \in \Xi^r$ such that $s(\uz)={\aa} $. The fact that  $(\Phi_i)_{i \in S}$ is a regular collection of functions  allows us to (locally) make a decomposition of the $\rho_r$ in terms of $\d {\aa} \times \d l_{\aa}$ (where $l_{\aa}$ is the induced volume measure on $\M_{\aa}^{[S]}$). This would enable us to refine the statement of Lemma 3.1 to the assertion that for $\mu(\aa)\d \aa$ a.e.-${\aa}$, it is true that for $\d l_{\aa}$-a.e. $\eta \in \M_{\aa}^{[S]}$, the events $\A(G_n,\eta)$ exhaust the support of $\P_\eta$.  We can then run the same argument as in Section \ref{singularity}, and conclude that for a.e. sub-manifold   $\M_{\aa}^{[S]}$ with identical values of the statistics corresponding to $(\Phi_i)_{i \in S}$, the Palm measures $\P_{\eta_1}$ and $\P_{\eta_2}$ are mutually singular for a.e. pair $(\eta_1,\eta_2) \in \M_{\aa}^{[S]} \times \M_{\aa}^{[S]}$.

\section{Acknowledgements}
The author is deeply grateful to the anonymous referees for their extremely valuable suggestions. This work was supported in part by the ARO grant W911NF-14-1-0094.

\end{document}